\documentclass[12pt,leqno]{amsart}
\usepackage[dvips]{graphics}
\usepackage{amssymb}

\numberwithin{equation}{section}

\theoremstyle{plain}
\newtheorem{thm}{Theorem}[section]
\newtheorem{lem}[thm]{Lemma}  
\newtheorem{prop}[thm]{Proposition}

\newtheorem{defn}[thm]{Definition}

\theoremstyle{remark}
\newtheorem{rem}{Remark}[section]
\newtheorem{ex}[rem]{Example}

\newcommand{\tref}[1]{Theorem~\ref{#1}}
\newcommand{\cref}[1]{Corollary~\ref{#1}}
\newcommand{\pref}[1]{Proposition~\ref{#1}}

\newcommand{\lref}[1]{Lemma~\ref{#1}}

\newcommand{\R}{\mathbb{R}}

\begin{document}
\title{Affine images of Riemannian manifolds}
\author{Alexander Lytchak}
\begin{abstract}
 We describe all affine maps from a Riemannian manifold to
 a metric space and all possible image spaces.
\end{abstract}


\thanks{The author was supported in part by the SFB 611 
{\it Singul\"are Ph\"anomene und Skalierung in mathematischen Modellen}, by the MPI for mathematics in Bonn
and by the Heisenberg grant from the DFG}

\maketitle
\renewcommand{\theequation}{\arabic{section}.\arabic{equation}}
\pagenumbering{arabic}

\section{Introduction}
A map $f:X\to Y$ between  metric spaces is called {\bf affine}
if it preserves the class of   linearly  parametrized minimizing geodesics. If $X$ is a geodesic metric space
then a continuous map $f$ is affine if and only if it sends midpoints
to midpoints. A map is called {\bf locally affine} (another word is {\bf totally geodesic}) if the restriction of the map to some neighborhood of any point
is affine.   Basic examples of (locally) affine maps are (locally) isometric embeddings, rescalings and projections to a factor in a direct product decomposition.

  Affine maps arise naturally in questions related to super-rigidity 
  (see the explanations and the literature list in \cite{Ohta}),
   in the study of Berwald spaces in Finsler geometry (\cite{Szabo}), in the study of product
  decompositions (\cite{deRham}) and in the study of isometric actions on non-positively curved spaces (\cite{Ballmann}).   
 It seems to me to be of independent interest, to understand  to what extent the geodesics determine the metric of a given space.

  We give a complete description of all affine maps 
  $f:M\to Y$ where $M$ is a smooth Riemannian manifold and $Y$ an arbitrary
  metric space. First,  an important special case of the main result:

  \begin{thm} \label{mainappl}
   Let $M$ be a connected complete smooth Riemannian manifold. There is a locally affine map
   $f:M\to Y$ to some metric space that is not a local homothety if and only if the universal covering 
   of $M$ is a product or a symmetric space of higher rank.
  \end{thm}

Here, we say that a map $f:X\to Y$ is  a \textbf{local homothety} if, for any point $x\in X$, there is some neighborhood 
$U$ of $x$ and a non-negative number $a$, such that for all $x_1,x_2 \in U$ the equality $d(f(x_1),f(x_2))=a\cdot d(x_1,x_2)$
holds true. If $X$ is connected then the number $a$ does not depend on the point $x$. 
Thus, up to rescaling by $a$, such  map locally is  an isometric embedding. 

For the formulation of the main theorem we will need two definitions.

 \begin{defn}
 A Riemannian submersion $f:M\to M_1$ between smooth Riemannian manifolds will be called
 \textbf{ flat} if the fibers are totally geodesic and the horizontal distribution is  integrable.
 \end{defn}
 
  A Riemannian submersion $f:M\to M_1$ between complete Riemannian manifolds is flat 
  if and only if the lift to the universal coverings $\tilde f :\tilde M\to \tilde M_1$ is the projection
  map of a direct product decomposition $\tilde M =\tilde M_1 \times M_2  \to \tilde M_1$. Another useful
  formulation is that a Riemannian submersion $f:M\to M_1$ is flat if and only if the tangent distribution of the 
  fibers is  a parallel distribution (\cite{vilms} or \cite{Sakai}).

 \begin{defn}
 Let $(M,g)$ be a Riemannian manifold. We will call a function
 $|\cdot | :TM \to \R$ a holonomy invariant Finsler structure
 if the restriction of $|\cdot |$ to each tangent space $T_p M$ is a 
(possibly not smooth or not strictly convex) norm and if $|v_1|=|v_2|$ for
any two vectors related by the holonomy along some piece-wise smooth curve.
For a holonomy invariant Finsler structure $|\cdot|$ we will call  
the identity map $Id: (M,g) \to (M,|\cdot |)$ an admissible change of metric.
 \end{defn}

 Note that an admissible change is a bi-Lipschitz map.  Clearly,  holonomy invariant
 Finsler structures are in one-to-one correspondence with norms $|\cdot|$ on a fixed tangent 
 space $T_p M$ that are invariant under the action of the holonomy group $Hol _p$. 
 We refer to Subsection \ref{description} for more information about    such norms.

 Now we can state the main result:
 
 \begin{thm} \label{main}
  Let $M$ be a complete Riemannian manifold. Then any locally affine map
  $f:M\to Y$ to a metric space $Y$ is a composition
  $f=f_i  \circ f_a \circ f_p$ of the following factors:
  \begin{enumerate}
  \item  The map $f_p :M\to M_1$ is flat Riemannian submersion onto a smooth complete Riemannian manifold $M_1$;
  \item The map $f_a :(M_1,g) \to (M_1,|\cdot |)$ is an admissible change of metric;
  \item The map $f_i$ is   a locally isometric embedding. 
 \end{enumerate}
 
  Moreover, any map $f$ of this kind is locally affine.
  \end{thm}

 This global theorem is a  consequence of the following
 local version of this result. To state it we recall  the following notation from \cite{Ohta}. 
 Let $f:M\to Y$ be a locally  affine map defined  on a Riemannian
 manifold $M$. For a tangent vector $v$ to $M$ we denote by $\gamma _v$ the
geodesic in the direction of $v$ and by $|v|^f$ the velocity of the image
geodesic $f\circ \gamma _v$.  Now the local result reads as follows:
 
 \begin{thm} \label{tech}
 Let $M$ be a smooth Riemannian manifold and let $f:M\to Y$ be a locally
 affine map. Then $|\cdot |^f$ is a continuous family of semi-norms on the
 tangent bundle of $M$ that is invariant under the parallel translation.
 \end{thm}

\begin{rem} It is important to stress that we consider all geodesics to be
parametrized affinely. If one parametrize all geodesics by the arclength
it is immediately clear that they determine the metric uniquely. If, on the 
other hand, one allows the geodesics to be parametrized arbitrary, then the question becomes much harder. For instance, if $X$ is flat and $Y$ a Finsler manifold this question is essentially Hilbert's fourth problem. If $X$ and $Y$ are Riemannian manifolds then such {\bf projective} mappings have recently been investigated  by Matveev (cf. \cite{Matveev}). In the setting of metric spaces it seems to be almost impossible to say something meaningful about such maps. 
\end{rem}

 The results of this paper generalize previous works  \cite{vilms}, \cite{Szabo}, \cite{Ohta} and \cite{Hitzel}. 
Vilms classifies in \cite{vilms} all affine maps between smooth Riemannian manifolds. Szabo extended  in \cite{Szabo}
this result 
to bijective maps between smooth Riemannian and smooth Finsler manifolds  providing a characterization of
all Berwald spaces.  On the other hand, Ohta (\cite{Ohta}) studied maps from a Riemannian manifold to  (locally uniquely geodesic)
metric spaces  and proved that such maps are compositions of affine maps to continuous Finsler manifolds and isometric
immersions. However, one cannot combine results of \cite{Ohta} and \cite{Szabo}, since the Finsler metric obtained is only 
continuous (even point-wise) and the Finsler geometric method used in \cite{Szabo} cannot work.

 Our proof is a self-contained combination and simplification of ideas used in \cite{Ohta} and \cite{Hitzel}. We hope,
 that our proof may lead to an understanding of affine maps on singular non-positively curved spaces.   

The paper is structured as follows.
In Section \ref{basics} we recall
 the proof of the semi-continuity of the the semi-Finsler structure $|\cdot| ^f$. In Section \ref{mainsec}, the core of the paper, we prove \tref{tech}.   In Section \ref{Prel} we recall  basics about invariant norms and Finsler structures.
Finally, in Section \ref{reduction} we prove \tref{mainappl} and \tref{main}.

\section{Basics} \label{basics}
\subsection{Setting} 
Let $M$ be a Riemannian manifold and let $f:M\to Y$ be a locally affine map to a metric space $Y$. We will denote
by $d$ the distance in $M$ and by $\bar d$ the distance in $Y$. In this and in the next section we are going to deal
with local questions only. Thus we may replace $M$ by a small strictly convex open neighborhood $U$ of a given point $x$.
If $U$ is sufficiently small then $f:U\to Y$ is affine. By the definition of $|\cdot |^f$, 
for any geodesic $\gamma _v :[a,b] \to V$ and all $s,t \in [a,b]$ we
 have $\bar d (f(\gamma _v (s)), f(\gamma _v (t))) = |s-t|\cdot |v|^f$.  From this we deduce 
 $|\lambda v|^f =|\lambda | \cdot |v|^f$  for all $\lambda \in \R$.

 In this section we are going to prove that 
$f$ is locally Lipschitz and that $|\cdot | ^f$  is a continuous family of semi-norms. With a minor additional assumption
it has been proved in \cite{Ohta}. For the convenience of the reader, we recall and slightly simplify Ohta's proof.  
 
\subsection{Lipschitz continuity}
We start with 
\begin{lem}
The map $f$ is continuous.
\end{lem}

\begin{proof}
We fix $p\in U$ and are going to prove the continuity at the point $p$.
In order to do so, it is enough to show that the homogeneous function $|\cdot |^f$ is bounded on the unit  sphere in $T^1 _p M  \subset T_p M$.  

  Choose $v_i \in T ^1 _p M$ to be the vertices of a regular  simplex $\Delta$. Set
  $x_i =\exp _p(\epsilon v_i)$ , for sufficiently small $\epsilon>0$. Let
  $A_0 (\epsilon )$ be the union of all $x_i$. Let $A_k (\epsilon )$ be the set of all points  that
  lie on a shortest geodesic between some pair of  points in $A_{k-1} (\epsilon ) $.
  The sets $K_{\epsilon} := \frac 1 {\epsilon} \exp _p ^{-1} (A_{n-1}  (\epsilon))$ converge in the Hausdorff topology
  to the boundary  $\partial \Delta$ of $\Delta$.  Thus $K_{\epsilon}$ does not contain the origin for small $\epsilon$.
  On the other hand, $K_{\epsilon}$ contains a small continuous perturbation of $\partial \Delta$. Thus $K_{\epsilon}$
  carries a non-trivial homology class of $H_{n-1} (T_p M \setminus \{ 0 \})$.  Therefore, $K_{\epsilon}$ intersects
  all rays starting from the origin.   
  Thus, for sufficiently small $\epsilon$, the compact set $A_{n-1} := A_{n-1} (\epsilon)$ does not contain $p$, but it intersects any geodesic  starting from $p$.
  
   Set $l=\max \{ \bar d(f(p),f(x_i)) \}$. By induction on $k$ and the triangle  inequality, $\bar d( f(p),f(x))\leq 2^k l$ for all $x\in A_k$.
   Since $A_{n-1}$ is compact, the number $\rho =d (p,A_{n-1})$ is positive.
 Thus for any $x\in A_{n-1}$ we get $\bar d(f(p),f(x))\leq C d(p,x)$ with
 $C=2^k \cdot l /\rho$. Since $A_{n-1}$ intersects any geodesic starting in
 $p$, we deduce that $|\cdot|^f$ is bounded by $C$ on $T_p ^1 M$.  
\end{proof}

\begin{lem}
The map $|\cdot|^f$ is continuous and $f$ is locally Lipschitz continuous. 
\end{lem}

\begin{proof}
 Choose  vectors $v_n \in T_{p_n} U$  converging to $v \in T_p U$.  Set $x_n=\exp (v_n); x=\exp (v)$. 
Then $x_n$ converges to $x$ and by continuity of $f$, the images $f(p_n)$
converge to $f(p)$ and the images $f(x_n)$ converge to $f(x)$.
Thus $|v_n|^f =\bar d(f(p_n),f(x_n))\to \bar d(f(p),f(x))=|v|^f$.

 The Lipschitz constant of $f$ at $p$ is bounded by the supremum of $|\cdot |^f$ on the unit sphere in $T_p U$.
 Thus, the continuity of $|\cdot |^f$ implies the Lipschitz continuity of $f$.    
\end{proof}

\subsection{Finsler structure} Now we claim:

\begin{lem}
The function $|\cdot|^f$ restricted to any tangent space $T_pU$ is a semi-norm.  
\end{lem}

\begin{proof}
 One can  follow the proof in (\cite{Ohta}). Another short proof  goes as follows.  The function $|\cdot |^f$ coincides with
 the metric differential defined by Kirhcheim (\cite{kirch})
 for any locally Lipschitz continuous map from a manifold to a metric space. Kirchheim proved that this metric differential is a semi-norm at almost all points. By continuity of $|\cdot |^f$ on $TU$, in our case this is a semi-norm at all points.
\end{proof}

\section{Main argument} \label{mainsec}
We are going to use the notations and assumptions from the last section.
In this section we are going to prove that the  semi-norms $|\cdot | ^f$ are  invariant 
 under parallel translations along arbitrary piecewise 
 $\mathcal C ^1$ curves in $U$. Approximating such a curve by a geodesic 
 polygon in the the $\mathcal C^1$-topology and using the fact that parallel
 transport behaves continuously with respect to such approximations, we deduce
 that it is enough to prove that the semi-norms are parallel along any geodesic polygon. Therefore, it is enough to prove that $|\cdot |^f$ is parallel along any geodesic $\gamma $ in $U$.
 
  Thus let us fix a geodesic $\gamma :[-a, a] \to U$ parametrized by the arclength. Set $p=\gamma (0)$ and  denote by 
  $P_t :T_p M \to T_ {\gamma (t)}M$ the parallel transport along $\gamma$.

 We call a vector $h \in T_x U$ \textbf{regular} if for all $v\in T_x U$ 
 $$ \lim _{t\to 0}\frac {|h+tv|^f +|h-tv|^f-2|h|^f} t =0$$

 A vector $h \in T_x U$ is regular if and only if  the semi-norm $|\cdot | ^f :T_x U \to \R$  is differentiable at $h$. 
  By the Theorem of Rademacher almost all vectors in $T_x M$ are  regular. Note that $h$ is regular if and only if 
  $rh$ is regular, for all  real non-zero $r$.

 We call a vector  $h\in T_p U$ \textbf{good} if, for almost
 all $t\in [-a,a] $,  the vector $h_t= P_t (h)$
 is a regular vector (in  $T_{\gamma (t)} U$). Applying Fubini's theorem
 we deduce that  almost all $h\in T_p U$ are good.  By continuity, it suffices to show that for any good vector $h$, the function $l(t)=|h_t|^f$ is constant.  Since $P_t$ and $|\cdot| ^f$ are positively homogeneous, we may assume after
 rescaling  that  $\exp (h_t) \in U$ exists  for all $t\in [-a,a]$.

  \begin{lem} The  function  $l(t)$ is Lipschitz continuous.     
   \end{lem}  
     
     \begin{proof}
     For all $t,s \in (-a,a)$,  by the triangle inequality $||h_t|^f-|h_s|^f|=
     |\bar d(f(\gamma (t)),f(\exp (h_t)) - \bar d (f(\gamma (s), f(\exp (h_s)))|  \leq 
     \bar d(f(\gamma (t), f(\gamma (s)) + \bar d (f(\exp (h_t), f(\exp (h_s))$.
     
     The parallel transport and the exponential map are  smooth and 
     $f$ is Lipschitz continuous, thus we can estimate $||h_t|^f-|h_s|^f|$ from above by $A\cdot |t-s|$,  for some $A>0$.
     \end{proof}

Since $l(t)$ is Lipschitz continuous, it is enough to prove
 that the derivative of $l(t)$ is $0$ almost everywhere. Thus, it suffices
 to prove that if $h_{t_0}$ is regular and if the derivative $ l '(t)$ exists at $t_0$, this derivative  must be zero.  To prove this claim, we may assume without loss of generality  (reparametrizing $\gamma$)  that $t_0 =0$ and that this derivative is non-negative.    
 Thus, we have reduced our task to proving the following claim:
 
 \begin{prop} \label{lastprop}
  Let $\gamma$ be a geodesic starting in $p$. Let $h$ be a regular  
  vector in  $T_p U$. Let $h_t$ be the parallel translates of $h$ along $\gamma$ and let the function $l(t) =|h_t|^f$ be differentiable at $0$.
  Then $l'(0) \leq 0$.
 \end{prop}

\begin{proof}
Again, we may assume that $\exp (h_t)$ exists for all $t$.
Assume $l'(0)= 2 \delta >0$. Then, for all small positive $t$, we have
$|h_t|^f\geq  |h|^f+ \delta t$.

Denote by $\eta _t $ the geodesic $\eta _t (r) =\exp (rh_t)$.
Set $x_{t,r}:= \eta _t (r)$ and $p=x_{0,0} = \gamma (0)$. Finally,
set $\mu _{t,r}:=\exp ^{-1} _p (x_{t,r}) \in T_pU$. The triangle inequality gives us
$$\bar d(f(x_{t,r}), f(x_{t,-r}) = 2r |h_t|^f \leq \bar d(f(p),f(x_{t,r})+\bar d(f(p)),f(x_{t,-r})=|\mu _{t,r}| ^f +
|\mu _{t,-r}|^f$$

Therefore
$$2r |h| ^f + 2r \delta t \leq |\mu _{t,r}| ^f + |\mu _{t,-r}|^f$$

 Denote by
$v_r \in T_pM$ the derivative $v_r := \frac d {dt} | _{t=0} \mu _{t,r}$.
The proof of the following  differential geometric lemma will be given below. 

\begin{lem}  \label{MainLemma}
In the above notations one has $\lim _{r\to 0} (||v_r -v_{-r}|| /r) =0$.
\end{lem}

 Given this lemma, it is easy to derive a contradiction:
 We find and fix a  small $r$ with $|v_r -v _{-r}| ^f \leq \delta r /4$.
 By definition, $\mu_ {t,\pm r} =   (\pm r)\cdot h +t \cdot v_{\pm r}  + o(t)$.
 Thus $|\mu _{t,-r}| ^f \leq |-rh + t \cdot v_r| ^f + \delta r t /2$, for small $t$.
 Hence $|\mu _{t,r}| ^f + |\mu _{t, -r}| ^f \leq |r\cdot h + t \cdot v_r| ^f+
 |rh -t \cdot v_r| ^f  +\delta r t $.
 
 By assumption, $h$ is regular. Thus $|rh+tv_r| ^f +|rh -t v_r|^f = 2 |rh|^f +o(t)$. For small $t$, this contradicts to 
 $|\mu _{t,r}|^f +|\mu _{t,-r}|^f \geq 2r + 2r \delta t$.
 \end{proof}

  It remains to provide:
  \begin{proof}[Proof of \lref{MainLemma}]
Let $\eta = \eta _0$. Denote by $Y$ the Jacobi field  $Y(r) = \frac d {dt} \eta_t(r) $ along $\eta$. Let $V_r$ be the Jacobi field  $V_r (s) := \frac d {dt} | _{t=0} exp _p (\frac {\mu _{t,r}} {r} s )$.

 By construction, $Y'(0)=0$, $V_r (0)=0$ and $V_r (r) =Y(r)$.
 Finally, $V_r ' (0) =\frac d {dt} \frac {\mu _{t,r}} {r}  = \frac {v_r} r$.
 Thus, we only need to prove that $|V_r '(0) +V_{-r} ' (0)|$ converges to
 $0$, as $r$ goes to $0$. Since the Jacobi field $(V_r+ V_{-r}) $ has value $0$ at the point $0$, it suffices to prove that $|(V_r + V_{-r} )(r)| =o(r)$.

  In fact we have $||Y(t)-Y(0)|| \leq C t^2$ and $||V_r (t) -t V_r ' (0)||  \leq C t^2$ for all $|t|\leq r$  and some $C$ independent of $r$. Together this  gives us $||V_{-r} (r) + V_{-r} (-r)|| \leq C r^2$ and $||Y(r)- Y(-r)|| \leq Cr^2$.
  Together this gives us $||V_{-r} (r) + V_r (r) || \leq C r^2$. This finishes the proof of \lref{MainLemma}, \pref{lastprop} and \tref{tech}.
    \end{proof}

\section{Remarks on invariant norms} \label{Prel} 
In this section we collect some elementary observation, probably well known to experts, for which we
could not find a reference.
\subsection{Invariant norms}
We recall a few definitions. A \textbf{semi-norm} on  a finite-dimensional vector space $V$ is a non-negative,
homogeneous,  convex function $q:V\to [0, \infty )$. It is called a \textbf{norm} if $q(v)>0$ for all $v\neq 0$.
It is called a Minkowski norm, if $q$ is smooth outside the origin and the Hessian $D^2_x q$ is positive
definite at all points $x\neq 0$. The second condition is equivalent to the requirement that 
$q$ can be expressed as $q=q_0+q_1$, where $q_0$ is a scalar product and $q_1$ is a   norm.    

For two norms  $q_1,q_2$ on $V$  there is some $L\geq 1$ with $q_1 \leq L \cdot q_2$ and $q_2 \leq  L \cdot q_1$.
The distance $|q_1-q_2|$ between $q_1$ and $q_2$ is defined to be the infimum of $\log (L)$ taken over the set of all such $L$.
The topology on the set of norms defined by this  distance function  coincides with the Hausdorff topology on the 
set of convex centrally-symmetric bodies.  It is well known that the set of Minkowski norms is dense in the set of all norms.

\begin{lem} \label{appr}
Let $V$ be a Euclidean space and let $G$ be a closed subgroup of the orthogonal group.
Then the set of $G$-invariant Minkowski norms is dense in the set of all $G$-invariant norms.
\end{lem}

\begin{proof}
Let $q$ be a $G$-invariant norm. Choose some small $\epsilon >0$.  We  find a Minkowski  norm $q_1$ with 
$|q_1 -q| \leq \epsilon$.  Let $q_2$ be the norm obtained from $q_1$ by the averaging procedure, i.e.,
$q_2 (x)= \int _G q_1 (gx) d \mu (g)$, where $\mu$ is the unit volume Haar measure on $G$.
 Then $q_2$ is  a $G$-invariant Minkowski norm and we still have $|q-q_2| \leq  \epsilon$, 
since $q$ is   $G$-invariant.  
\end{proof}

\begin{lem}
Let $V, G$ be as above. If $G$ acts transitively on the unit sphere then the only $G$-invariant norms are multiples
of the given Euclidean norm. If $G$ does not act transitively on the unit sphere then there are non-Euclidean 
$G$-invariant Minkowski norms.
\end{lem}

\begin{proof}
The first statement is clear. To prove the second statement we proceed as follows. If $G$ acts irreducibly, consider
any orbit $Gp$ and let $C$ be its convex hull. If $G$ acts reducibly, we consider  orthogonal $G$-invariant
subspaces $V_1,V_2 \subset V$ with $V_1 \oplus V_2 =V$, and let $C$ be the convex hull of the union of the unit spheres $S_1 \subset V_1$ and
$S_2 \subset V_2$. In both cases we obtain a $G$-invariant convex body that is not strictly convex. Its symmetrization
around the origin and  the norm $q$ defined by the symmetrized body is still $G$-invariant and not strictly convex.
Due to  the previous lemma, we find a $G$-invariant Minkowski norm $q_n$ arbitrary close to $q$. Since $q$ is non-Euclidean,
$q_n$ is non-Euclidean as well, at least for large $n$. 
\end{proof}

\subsection{Holonomy group} For simplicity,
we will  restrict ourselves to the complete case. Thus let $(M,g)$ be a complete smooth 
Riemannian manifold. Let $p\in M$ be a point and let $H$ and $H^0$ denote the holonomy group $Hol _p$ and its identity component
respectively. The tangent space $V=T_p M$ splits under the action of $H^0$ in a  uniquely defined way as
$V=V_0 \oplus V_1 \oplus ...\oplus V_i$, where the action of $H$ on $V_0$ is trivial and the action on $V_i, i\geq 1$ is 
irreducible and not trivial. Moreover, this decomposition determines the unique  direct product  decomposition of the universal
covering $\tilde M$ of $M$ (\cite{Sakai}).
 By the uniqueness of this decomposition up to permutation, the action of $H$ preserves
this decomposition of $V$, possibly up to a   permutation of summands. In particular, if $H^0$ acts reducibly then the closure
$\bar H$ of $H$ acts non-transitively on the unit sphere.

 If the action of  $H^0$ on the unit sphere is irreducible, then either it acts transitively on the unit sphere
 or the space $M$ is a locally symmetric space of rank at least two, with irreducible universal covering, due to
 the theorem of Berger-Simons (\cite{simons}). In this case, each
 element of $H$ acts as the differential of a local isometry, thus it preserves the type of a vector in its Weyl chamber.
 Hence, the action of $\bar H$ on the unit sphere is non-transitive as well.

 We conclude:
 \begin{lem}  \label{berger}
 Let $M$ be a complete Riemannian manifold. Let $p\in M$ be a point, $V=T_p M$ and let $H$ be the holonomy group at
 the point $p$. Then the following are equivalent:
 \begin{enumerate}
 \item  $H$ does not act transitively on the unit sphere $V^1$ of $V$;
 \item The closure of $H$  does not act transitively on $V^1$;
 \item There are non-Euclidean  Minkowski norms on $V$ invariant under $H$;
 \item The universal covering of $M$ is either a direct product or a symmetric space of higher rank.
 \end{enumerate}
 \end{lem} 

\subsection{Description of invariant norms} \label{description}
In this subsection we are going to provide a description of (Minkowski) norms invariant
under the action of a continuous holonomy group.   Since these statements are not used in the rest of the paper we will
omit some details. 

 Thus let $M$ be a Riemannian manifold and assume  that the holonomy group $H=Hol_p$ is connected (this happens, for instance,
 if $M$ is complete and simply connected). Then (due to the theorems of Berger-Simons an de Rham)
 the action of $H$ on $V=T_p M$ is polar, i.e., there is a linear subspace $\Sigma$ of $V$ that intersects all
orbits of $H$ with all intersection being orthogonal. The stabilizer $N(\Sigma )$ acts on $\Sigma $ as a finite Coxeter
group $W$.  Then each $H$-invariant norm $q$ on $V$ restricts to a $W$-invariant norm $\bar q$  on $\Sigma$. We have:

\begin{prop}
 The restriction $q\to \bar q$ is a bijection between the set of $H$-invariant norms on $V$ and $W$-invariant
 norms on $\Sigma$. Moreover,  $q$ is a Minkowski norm if and only if $\bar q$ is a Minkowski norm. 
\end{prop}

\begin{proof}
 If $q$ is a (Minkowski) norm then so is its restriction to any subspace, in particular $\bar q$. On the other hand,
 each $W$-invariant function $\bar q$ extends to a unique $H$-invariant function $q$ on $V$. Now,  if $\bar q$
 is a norm then so is $q$, due to \cite{chi}, p. 107.  Moreover,  if one can represent $\bar q$ as a sum
 $\bar q_0 + \bar q_1$ of a scalar product and a norm, then, averaging, we may assume that $\bar q_0$ and $\bar q_1$
 are $W$-invariant as well. Thus they extend to an $H$-invariant scalar product and norm $q_0$ and $q_1$
with $q=q_0+q_1$.   Finally, an $H$-invariant function is smooth if and only if its restriction to $W$ is smooth
(\cite{dadok}). Smoothing $\bar q$ at the origin, 
we deduce, that $\bar q$ is smooth in $\Sigma \setminus \{ 0\} $ if and only
if $q$ is smooth in $V\setminus \{ 0\}$.  Thus $q$ is a Minkowski norm if and only if $\bar q$ is. 
\end{proof}

We illustrate this result by two examples:

\begin{ex}
Let $M$ is a symmetric space and  let $F$ be a maximal flat through a point $p$. Let $W$ be the group of all isometries
of $M$ that fix $p$ and leaves $F$ invariant. Then the set of parallel (smooth)  Finsler structures is in one-to-one correspondence with the set of all $W$-invariant (Minkowski) norms on $T_p M$.   
\end{ex}

\begin{ex}
Let $M=M_1\times M_2$ be the product of two irreducible not locally symmetric spaces. For $p\in M$, choose
any unit vector $e_1$ in $T_p M_1$ and $e_2$ in $T_p M_2$.  Then one can choose $\Sigma$ to be the plane
generated by $e_1$ and $e_2$.  The group $W$ has  $4$ elements and is generated by the the reflections
$e_i\to - e_i$. 
\end{ex}

\section{The conclusion} \label{reduction}
\subsection{The if part}
We are going to prove \tref{main}. First we would like to see that each map as described in the theorem is locally affine.

A composition of continuous locally affine map is locally affine. A locally isometric embedding is locally affine.
A flat submersion is locally given by a projection onto a direct factor, thus it is locally affine as well. Therefore, it suffices to prove that an admissible change of a metric is locally affine. We are going to reduce the statement to the smooth
case, where the result is known (cf. \cite{Szabo}).  

 Thus let $(M,g)$  be a smooth complete Riemannian manifold and   let $|\cdot |$ be a holonomy invariant Finsler structure.
 Fix a point $x\in M$. Aproximate the   norm  $|\cdot|$ on $T_p M$ by Minkowski norms   $|\cdot |^n$  norms on $T_p M$  invariant under the holonomy group  $Hol _p$ (\lref{appr}).  Extend each norm $|\cdot|^n$ to a smooth holonomy invariant
 Finsler structure on $M$. Due to \cite{Szabo}, the identity map $Id: (M,g) \to (M, |\cdot |^n )$ is affine.
 Since minimizing
 geodesics in the smooth Finsler manifold $(M,|\cdot|^n)$ are locally unique,  all $(M,|\cdot |^n)$ geodesics are 
 images of $(M,g)$ geodesics (i.e., the map $Id: (M,|\cdot |^n ) \to (M,g)$ is affine as well).
 We find some $L$ such that the identity $Id: (M,g) \to (M,|\cdot|^n)$ is L-bilipschitz, for all $n$.
We find some strictly convex ball $V$ of radius $r$ around $x$, such that any geodesic between points of $V$ that
is not contained in $V$ has length at least $4Lr$. It follows, that any $(M,g)$-geodesic  contained in $V$
is a minimizing geodesic in $(M,|\cdot|^n)$, for all $n$. Since minimizing (!) $(M,|\cdot| ^n)$ geodesics converge
to minimizing $(M,|\cdot |)$ geodesics, we deduce that $Id:(V,g) \to (V,|\cdot| )$ sends minimizing geodesics
to minimizing geodesics. The linear parametrization of the geodesics is clear by approximation (or by the fact that
geodesics are auto-parallel).

   \subsection{Decomposition of $f$} 
Let $M$ be a complete Riemannian manifold and let $f:M\to Y$ be a locally affine map.   Due to \tref{tech}, the function
$|\cdot |^f$ defines a parallel family of semi-norms. For each $q\in M$, we set 
$V_q = \{ v\in T_q M | |v|^f =0 \}$.  Since the semi-norms are parallel this is a parallel distribution. 

 Assume first that $M$ is simply connected. Then by the theorem of de Rham this distribution $q\to V_q$ is the vertical distribution
 of a projection $f_p:M=M_1\times M_2 \to M_1$ onto a direct factor. The map $f$ factors through $f_p$ as $f=\hat f \circ f_p$.
 Moreover, $\hat f$ coincides with the restriction of $f$  to any horizontal slice $M_1 \times \{ x_2\}$. Thus $\hat f$
 is again affine, and the semi-norm $|\cdot|^{\hat f}$ coincides with the restriction of $|\cdot|^f$. By the definition of
 $V_q$, this semi-norm is a norm.  Thus on $M_1$ we have an admissible change of metric $(M_1,g) \to (M_1,|\cdot |^ f  )$
 and the induced map $f_i :(M_1, |\cdot |^f) \to Y$ is by definition of $\cdot| ^f$ a local isometry.

Assume now that $M$ is not simply connected. Consider the lift $\tilde f$ of $f$ to the universal covering $\tilde M$.
Consider a horizontal slice $M'_1$ in $\tilde M$. The restriction of $\tilde f$ to this slice is locally bi-Lipschitz. 
On the other hand this restriction factors through the canonical projection of the slice $M' _1$ to $M$. This implies that
the leaves of the integral distribution $V_q$ in $M$  are closed submanifolds that define a Riemannian  submersion $f_p:M\to M_1$
for some manifold $ M_1$ covered by $M' _1$. It is clear that the Riemannian submersion $f_p$ is flat. The rest follows from the simply connected case.   

\subsection{Proof of \tref{mainappl}}
To prove \tref{mainappl}, first let $M$ be  locally irreducible and not locally symmetric of higher rank. Let $f:M\to Y$ be a locally affine map.  Due to \tref{tech} and \lref{berger},  at each point $p\in M$, the semi-norm $|\cdot |^f$ on $T_p M$  is either constant $0$ or a positive multiple
of the metric $g$, the multiple not depending on the point. Thus the map $f$ either sends the whole manifold $M$ to a point
or rescales distances in any small neighborhood of any given point by the constant number $a$. 
 
  If, on the other hand, $M$ is either locally reducible or locally symmetric of higher rank, then there is an admissible change
  of metric to a non-Euclidean (smooth) Finsler structure. This change is certainly not a local homothety,  but it is a
  locally affine map by \tref{main}.


\end{document}